\documentclass{etds}

\numberwithin{equation}{section}

\font\Bbb=msbm10 scaled 1000
\def\sym{\fam\comfam\com}
\font\tensym=msbm10 at 12pt \font\sevensym=msbm7
\font\fivesym=msbm5 
\newfam\symfam
\textfont\symfam=\tensym \scriptfont\symfam=\sevensym
\scriptscriptfont\symfam=\fivesym
\def\sym{\fam\symfam\relax}
\def\tower{
\setlength{\unitlength}{1mm}
\begin{picture}(0,7)
\put (0,0) {\framebox (50,20)} \put (3,-4) {\makebox (3,-4)[br]
{$B_{k}$}} \put (16,-6) {\makebox (16,-6)[b]{$\underbrace {
   \quad  \quad   \quad  \quad \quad  \quad  \quad  \quad  \quad  \
  \quad  \quad \quad  \quad
  }_{B_{k-1}}$}}
\put (17,-16) {\makebox (17,-16)[b]{$\overline {
      \quad  \quad  \quad \quad ~ ~ p_{k-1}
\quad  \quad   \quad  \quad \quad   \quad ~ ~
  }$}}
\multiput (5,0)(5,0){10}{\line (0,20){20}} \put (0,4) {\line
(4,0){50}} \put (3,0){\vector (0,1){4}} \thicklines \put (0,0)
{\line (5,0){5}} \multiput (0,21)(0,1){5}{\line (5,0){5}} \put
(-3,27){\makebox (-2,27)[bl] {$a_1^{(k-1)}$}} \put
(5,21){\line(5,0){5}} \put (5,22){\line(5,0){5}} \put
(8,24){\makebox (8,24)[bl] {$a_2^{(k-1)}$}} \put (10,21){\makebox
(10,21)[b] {$\cdots \cdots$}} \multiput (25,21)(0,1){6}{\line
(25,0){5}} \put (23,28){\makebox (23,28)[bl] {$a_j^{(k-1)}$}} \put
(26,21){\makebox (26,21)[b] {$\cdots \cdots$}} \multiput
(45,21)(0,1){3}{\line (45,0){5}} \put (45,25){\makebox (45,25)[bl]
{$a_{p_{k-1}}^ {(k-1)}$}} \put (-40,-25){\centerline{Figure 1 :
$k^{\hbox {th}}$--tower.}}
\end{picture}}

\def\proof #1 {\noindent {\bf Proof #1 : }}  

\def\build#1_#2^#3{\mathrel{\mathop{\kern 0pt#1}\limits_{#2}^{#3}}}
\def\tend#1#2{\build\hbox to 12mm{\rightarrowfill}_{#1\rightarrow #2}^{}}
\def\netendpas#1#2{\build\hbox to 12mm{$\not \longrightarrow$}_{#1 \rightarrow
#2}^{}}

\def\tendn{\tend{n} {\infty}}

\def \converge#1#2#3{\build\hbox to 15mm {\rightarrowfill}_{#1\rightarrow #2}^
{\hbox{\scriptsize #3}}}

\def\N{{\sym N}}
\def\Z{{\sym Z}}

\def \T{{\sym T}}
\def \P{{\sym P}}
\def \E{{\sym E}}

\begin{document}
\ETDS{1}{1}{?}{2006}

\runningheads{E.\ H. \ El Abdalaoui}{rank one transformations with
singular spectrum}

\title{A new class of rank one transformations with singular
spectrum $^{*}$}

\author{EL \ HOUCEIN EL \ ABDALAOUI}

\address{Universit\'e de Rouen-Math\'ematiques \\
  Labo. de Maths Raphael SALEM  UMR 60 85 CNRS\\
Avenue de l'Universit\'e, BP.12 \\
76801 Saint Etienne du Rouvray - France . \\
\email{elhoucein.elabdalaoui@univ-rouen.fr}}

\recd{$31$ July $2006$}

$^*$ Dedicated to Professor J. De Sam Lazaro\\

\begin{abstract}
We introduce a new tool to study the spectral type of rank one
transformations using the method of central limit theorem for
trigonometric sums. We get some new applications.
\end{abstract}

\section{Introduction}
~~~~~The purpose of this paper is to bring a new tool in the
study of the spectral type of rank one transformations. Rank one
transformations have simple spectrum and in \cite{Ornstein} D.S.
Ornstein, using a random procedure, produced a family of mixing
rank one transformations. It follows that the Ornstein's class of
transformations may possibly contain a candidate for Banach's well-known
problem whether there exists a dynamical system $(\Omega,{\mathcal
{A}},\mu,T)$ with simple Lebesgue spectrum. But, in 1993, J.
Bourgain in \cite{Bourgain} proved that almost surely Ornstein's
transformations have singular spectrum. Subsequently, using the
same method, I. Klemes \cite{Klemes1} and I. Klemes \& K. Reinhold
  \cite{Klemes2} obtain that mixing 
staircase transformations of Adams \cite{Adams1} and Adams \&
Friedman \cite{Adams2} have singular spectrum. They conjecture
that all rank one transformations have singular spectrum.\par

Here we shall exhibit a new class of rank one transformations with
singular spectrum. Our assumption include some new class of
Ornstein transformations and a class of Creutz-Silva rank one
transformations \cite{Creutz-silva}. Our proof is based on
techniques introduced by J. Bourgain
\cite{Bourgain} in the context of rank one transformations and
developed by Klemes \cite{Klemes1}, Klemes-Rienhold
\cite{Klemes2}, Dooley-Eigen \cite{Eigen}, together with some ideas from the proof of 
the central limit theorem for trignometric sums. The fundamental key, as noted by Klemes
\cite{Klemes1}, is the estimation of the $L^1$-norm of a certain
 trigonometric polynomial ${{(|P_m|^2-1)}}$
. We shall use the method of central limit theorem for trignometric sums to produce an
estimate of this $L^1$-norm.

\section*{2. Rank One Transformation by Construction}

Using the cutting and stacking method described in [Fr1], [Fr2],
one can construct inductively a family of measure preserving
transformations, called rank one transformations, as follows
\vskip 0.1cm Let $B_0$ be the unit interval equipped with the
Lebesgue measure. At stage one we divide $B_0$ into $p_0$ equal
parts, add spacers and form a stack of height $h_{1}$ in the usual
fashion. At the $k^{th}$ stage we divide the stack obtained at the
$(k-1)^{th}$ stage into $p_{k-1}$ equal columns, add spacers and
obtain a new stack of height $h_{k}$. If during the $k^{th}$ stage
of our construction  the number of spacers put above the $j^{th}$
column of the $(k-1)^{th}$ stack is $a^{(k-1)}_{j}$, $ 0 \leq
a^{(k-1)}_{j} < \infty$,  $1\leq j \leq p_{k}$, then we have

$$h_{k} = p_{k-1}h_{k-1} +  \sum_{j=1}^{p_{k-1}}a_{j}^{(k-1)}.$$
\vskip 3 cm \hskip 3.5cm
  \tower
\vskip 3.0cm
\noindent{}Proceeding in this way we get a rank one transformation
$T$ on a certain measure space $(X,{\mathcal B} ,\nu)$ which may
be finite or
$\sigma-$finite depending on the number of spacers added. \\
\noindent{} The construction of a rank one transformation thus
needs two parameters, $(p_k)_{k=0}^\infty$ (parameter of cutting
and stacking), and $((a_j^{(k)})_{j=1}^{p_k})_{k=0}^\infty$
(parameter of spacers). Put

$$T \stackrel {def}= T_{(p_k, (a_j^{(k)})_{j=1}^{p_k})_{k=0}^\infty}$$

\noindent In \cite{Nadkarni1} and \cite{Klemes2} it is proved that
the spectral type of this transformation is given (upto possibly some discrete measure) by

\begin{eqnarray}\label{eqn:type1}
d\sigma  =W^{*}\lim \prod_{k=1}^n\left| P_k\right| ^2d\lambda,
\end{eqnarray}
\noindent{}where
\begin{eqnarray*}
&&P_k(z)=\frac 1{\sqrt{p_k}}\left(
\sum_{j=0}^{p_k-1}z^{-(jh_k+\sum_{i=1}^ja_i^{(k)})}\right),\nonumber  \\
\nonumber
\end{eqnarray*}
\noindent{}$\lambda$ denotes the normalized Lebesgue measure on the
circle group $\T$ and $W^{*} \lim$ denotes weak*limit in the space of
bounded Borel measures on ${\T}$.\\

\noindent{}The principal result of this paper is the following:\\

\noindent {\bf Theorem 2.1. }{\it Let $T=T_{(p_k,
(a_j^{(k)})_{j=1}^{p_k})_{k=0}^\infty}$ be a rank one
transformation such that,

\begin{eqnarray}{\label{eq:eqt}}
&&(i)~~a_{j+1}^{(k)}\geq 2s_k(j),j=0,\cdots, p_k-1,\nonumber\\
&& {\rm {~~with~~~~}} s_k(j)=a_1^{(k)}+\ldots +a_j^{(k)},
s_k(0)=0,\nonumber\\
&&(ii)~~\frac{s_k{(p_k)}}{h_k}<\frac12 \nonumber
\end{eqnarray} then the spectrum of $T$ is singular.}
\\

\noindent{}We remark that the
spectrum of rank one transformation is always singular if the
cutting parameter $p_k$ is bounded. In fact, Klemes-Reinhold
prove that if $\displaystyle \sum_{k=0}^{\infty}
\frac{1}{{p_k}^2}=\infty$ then the associated rank one
transformation has singular spectrum. Henceforth, we assume that the series
$\displaystyle \sum_{k=0}^{\infty} \frac{1}{{p_k}^2}$
converges.\\
We point out also that the condition $(ii)$ of
theorem 2.1 holds in the case of rank one transformations satisfying a
restricted growth condition  provided that $\displaystyle \min_{1 \leq j \leq  p_k}(a_j^{(k)})=0$. Following Creutz-Silva \cite{Creutz-silva},
we say that a rank one transformation $T=T_{(p_k, (a_j^{(k)})_{j=1}^{p_k})_{k=0}^\infty}$ has
restricted growth if
\[
\displaystyle \frac{\displaystyle \sum_{i=1}^{p_k} \left( a_j^{(k)}- \displaystyle \min_{1 \leq j \leq  p_k}(a_j^{(k)})\right) }{h_k+\displaystyle \min_{1 \leq j \leq  p_k}(a_j^{(k)}) }\tend{k}{+\infty}0,
\]
\noindent{}The proof of our main result is based on the method of
J. Bourgain in \cite{Bourgain}. We shall recall the main ideas of
this
method.\\

\noindent {\bf Proposition 2.2. }{\it \ The following are
equivalent  }

\begin{enumerate}
\item [(i)]  {\it $\sigma \perp \lambda$ }

\item[(ii)]  {\it $\inf \{\displaystyle \displaystyle \int \displaystyle \int
\prod_{l=1}^k\left| {P_{n_l}(z)}\right| d\lambda,~k\in
{\N},~n_1<n_2<\ldots <n_k\}=0.$ }
\end{enumerate}

\noindent{}Now fix some subsequence ${\mathcal {N}}=\left\{
n_1<n_2<\ldots <n_k\right\}$ ,\linebreak[0]$\ k \in {\N}$, $m>n_k$
and put
\[
Q\left(z\right) =\prod_{i=1}^k\left |{P_{n_i}(z)}\right |.
\]

\noindent{}One can show, using the same arguments as in
\cite{elabdal1}, the following lemma.\\

\noindent {\bf Lemma 2.3. }{\it
\begin{eqnarray*}
  &&\displaystyle \int  Q(z) \left| P_m(z)\right|d\lambda(z)
\leq \\
  &&\frac 12\left( \displaystyle \int Q d\lambda +\displaystyle \int
Q(z) \left| P_m(z)\right| ^2 d\lambda(z)\right) -\frac 18\left(
\displaystyle \int  Q \left| \left| P_m(z)\right| ^2-1\right|
d\lambda(z)\right) ^2. \end{eqnarray*}}\\

\noindent {\bf Proposition 2.4. }{\it $\displaystyle \lim_{m
\rightarrow
  \infty}\displaystyle \int  Q \left|
P_m(z)\right|^2d\lambda(z) {\Bbb {=}} \displaystyle \displaystyle
\int  Q $ $d\lambda(z)  .$}\\
\\
\begin{proof}
\noindent{}The
sequence of probability measures $\left| {P_m(z)}%
\right|^2d\lambda(z)$ converges weakly to the Lebsegue measure.
   \end{proof}
   \\

\noindent We have also the following proposition: \\

\noindent {\bf Proposition 2.5. }{\it There exist a subsequence of
the sequence $\left (\left |\left| P_m(z)\right|-1 \right|\right)$
which converge weakly to some non-negative function $\phi$ which
satisfies  $ \phi \leq $2, almost surely with respect to the Lebesgue
measure.}
\\

\begin{proof}
\noindent{}The
sequence $\left |\left| {P_m(z)}%
\right|-1 \right |$ is bounded in $L^2$. It follows that there
exist a subsequence which converges weakly to some non-negative
$L^2$ function $\phi$. Let $\omega$ be a non-negative continuous
function, then we have

\begin{eqnarray*}
\int \omega\left |\left| {P_m(z)} \right|-1 \right |d\lambda(z)
&&\leq \int \omega\left| {P_m(z)} \right| d\lambda(z)+\int \omega
d\lambda(z)\\
&&\leq {(\int \omega d\lambda(z))}^{\frac12} {(\int \omega \left |
{P_m(z)} \right|^2d\lambda(z))}^{\frac12}+\int \omega d\lambda(z).
\end{eqnarray*}
\noindent{}Hence
\[
\int \omega \phi d\lambda \leq 2 \int \omega d\lambda.
\]

Since this holds for all non-negative continuous $\omega$, we have $\phi \leq 2$ a.e. 
\end{proof}
\\

\noindent Put
\[
\alpha = \phi~~ d\lambda.
\]
\noindent We shall prove the following:\\

\noindent {\bf Proposition 2.6. }{$\alpha \bot \sigma$}.\\

\noindent For the proof of the proposition 2.6 we need the
following classical lemma \cite{Kilmer-Seaki}.\\

\noindent {\bf Lemma 2.7. }{\it Let $\rho, \tau$ be two
nonnegative finite measure on a measurable space $X$. Then the
following properties are equivalent :
\begin{enumerate}
\item[(i)]  {\it $\rho \perp \tau$ }

\item[(ii)]  {\it Given $\varepsilon>0$, there exists a
nonnegative measurable function $f$ on $X$ such that $ f >0,$
$\tau-a.e.$ and such that
\[
\left(\int f d\rho) \right) \left (\int \frac{d\tau}{f} \right)
<\varepsilon.
\] }
\end{enumerate}
}

\begin{proof}{of Proposition 2.6}
Let $\displaystyle \beta_1=\inf\{\int Q
d\alpha,~:~Q=\prod_{i=1}^k\left | {P_{n_i}(z)}\right|, k \in \N,
n_1<n_2<\cdots<n_k\}$ and $\displaystyle \beta_2=\inf\{\int Q
d\lambda,~:~Q=\prod_{i=1}^k\left | {P_{n_i}(z)}\right|, k \in \N,
n_1<n_2<\cdots<n_k\}$. Then, for any $Q$, we have

\begin{eqnarray*}
    \int Q
d\alpha \geq \beta_1,~~~ \liminf \int Q |P_m| d\lambda \geq
\beta_2
\end{eqnarray*}

\noindent{} Combine lemma 2.3, proposition 2.4 and 
proposition 2.5 and take the inf over all $Q$ to obtain:
\begin{eqnarray*}
\beta_2 \leq \beta_2-\frac18\beta_1^2
\end{eqnarray*}
\noindent It follows that
\[
\beta_1=0.
\]
\noindent{}We claim: $\beta_1=0$ implies $\displaystyle \int
\prod_{j=1}^{k}|P_j|d\alpha \tend{k}{\infty}0$. In fact, let
${\mathcal{N}}=\{n_1<n_2<\cdots<n_k\}\subset \N^*$ and $N>n_k$.
Then, by Cauchy-Schwarz inequality, we have

\begin{eqnarray}{\label{eq:Bo1}}
\int \prod_{j=1}^{N}|P_j|d\alpha & =&  \int \sqrt{\prod_{j \in
\mathcal{N}}|P_j|} \sqrt{\prod_{j \in \mathcal{N}}|P_j|} \prod_{j
\not \in \mathcal{N}}|P_j| d\alpha \nonumber\\
&\leq& {\left(\int \prod_{j \in \mathcal{N}}|P_j|d\alpha\right
)}^{\frac12} {\left(\int \prod_{j \in \mathcal{N}} |P_j| d\alpha
\prod_{j \not \in \mathcal{N}}|P_j|^2 d\alpha\right )}^{\frac12}
\end{eqnarray}

\noindent{}But

\begin{eqnarray}{\label{eq:B2}} \int \prod_{j \in \mathcal{N}}
|P_j|  \prod_{j \not \in \mathcal{N}}|P_j|^2 d\alpha &\leq & 2
\int \prod_{j \in \mathcal{N}} |P_j|  \prod_{j \not \in
\mathcal{N}}|P_j|^2 d\lambda \nonumber \\
&\leq&  2\int \prod_{j=1}^{N} |P_j|  \prod_{j \not \in
\mathcal{N}}|P_j| d\lambda \nonumber\\
  &\leq& 2{\left(\int \prod_{j=1}^{N} |P_j|^2 d\lambda \right)}^{\frac12}
{\left(\int\prod_{j \not \in \mathcal{N}}|P_j|^2
d\lambda\right)}^{\frac12} \nonumber \\
&\leq& 2
\end{eqnarray}

\noindent{}Combine (\ref{eq:Bo1}) and (\ref{eq:B2}) to get the
claim. The proof of the proposition follows from lemma 2.7.
\end{proof}

\section*{3. Estimation of $L_1$-Norm of ($\left |P_m(z)\right |-1$)
and the Central Limit Theorem.}

In this section we assume that for $m$
sufficiently large
\begin{eqnarray}{\label{eq:rg}}
  &&(i)~~a_{j+1}^{(m)}\geq 2s_m(j),j=0,\cdots, p_m-1,\nonumber
  \\ &&{\rm {~~with~~~~}} s_m(j)=a_1^{(m)}+\ldots +a_j^{(m)},
s_m(0)=0,\nonumber\\
&&(ii)~~\frac{s_m{(p_m)}}{h_m}<\frac12 \nonumber
\end{eqnarray}

\noindent{}Under the above assumptions, we shall proof the following:\\

\noindent {\bf Proposition 3.1. }{ $\alpha \geq K \lambda$, for
some positive
constant $K$}.\\

\noindent{} The proof of the proposition is based on an
estimate of $\displaystyle \int_A ||P_m|-1|d\lambda$,   where
$A$ is a Borel set with $\lambda(A)>0$. More precisely we shall
study the stochastic behavior of the sequence $ |P_m|$. For that
our principal strategy is based on the method of the central
limit theorem for trigonometric sums. A nice account can be
found in \cite{Kac}. It is well-known that
Hadamard lacunary trigonometric seires satisfies the central limit
theorem \cite{Zygmund}. The central limit theorem for
trigonometric sums has been studied by many authors, Zygmund and
Salem \cite{Zygmund}, Erd\"os \cite{Erdos},
  J.P. Kahane \cite{Kahane}, Berkers \cite{Berkes}, Murai \cite{Murai},
Takahashi \cite{Takahashi}, Fukuyama and  Takahashi
\cite{Fukuyama}, and many others. The same method is used to
study the asymptotic behavior of Riesz-Raikov sums \cite{petit}.\\

\noindent Here, we shall prove the following:\\

\noindent {\bf Proposition 3.2. } {If $(i)$ and $(ii)$ hold then
for any Borel subset $A$ of $\T$ with $\lambda(A)>0$, the distribution of the
sequence of random variables
$\frac{\sqrt{2}}{\sqrt{p_m}}\sum_{j=0}^{p_m-1} \cos((jh_m+s_m(j))t)$
converges to the Gauss distribution. That is

\begin{eqnarray}{\label{eq:eq6}}
\frac1{\lambda(A)}\left \{t \in A~~:~~\frac{\sqrt{2}}{\sqrt{p_m}}
\sum_{j=0}^{p_m-1} \cos((jh_m+s_m(j))t) \leq x \right \}\nonumber \\
\tend{m}{\infty}\frac1{\sqrt{2\pi}}\int_{-\infty}^{x}
e^{-\frac12t^2}dt\stackrel{\rm {def}}{=}{\mathcal {N}}\left( \left] -\infty,x \right] \right) .
\end{eqnarray}
}\\
The proof of proposition 3.2 is based on the idea of the proof of
martingale central limit theorem due to McLeish \cite{Mcleish}.
The main ingredient is the following .\\

\noindent {\bf Lemma 3.3. } {For $n \geq 1$, let $U_n,T_n$ be
random variables such that
\begin{enumerate}
\item  $U_n\longrightarrow a$ in probability.
\item $\{T_n\}$ is uniformly integrable.
\item $\{|T_nU_n|\}$ is uniformly integrable.
\item $\E(T_n) \longrightarrow 1$.
\end{enumerate}
Then $\E(T_nU_n)\longrightarrow a$.}\\

\begin{proof}
Write $T_nU_n=T_n(U_n-a)+aT_n$. As $\E(T_n) \longrightarrow
1$, we need to show that $\E(T_n(U_n-a))\longrightarrow 0$.\\
Since $\{T_n\}$ is uniformly integrable, we have
$T_n(U_n-a)\longrightarrow 0$ in probability. Also, both  $T_nU_n$
and $aT_n$ are uniformly integrable, and so the combination
$T_n(U_n-a)$ is uniformly integrable. Hence,
$\E(T_n(U_n-a))\longrightarrow 0$.
\end{proof}\\

\noindent{}Let us recall the following expansion
\[
\exp(ix)=(1+ix)\exp(-\frac{x^2}2+r(x)),
\]
\noindent{}where $|r(x)| \leq |x|^3$, for real $x$.\\

\noindent {\bf Theorem 3.4. }{Let $\{X_{nj}~~:~~ 1 \leq j \leq
k_n, n\geq 1\}$ be a triangular array of random variables.
$\displaystyle S_n=\sum_{j=1}^{k_n}X_{nj}$, $\displaystyle T_n=\prod_{j=1}^{k_n}(1+itX_{nj})$,
and $\displaystyle U_n=\exp \left (-\frac{t^2}2 \sum_j
X_{nj}^2+\sum_{j}r(tX_{nj})\right )$. Suppose that
\begin{enumerate}
\item $\{T_n\}$ is uniformly integrable.
\item $\E(T_n) \longrightarrow 1$.
\item $\displaystyle \sum_j X_{nj}^2 \longrightarrow 1$ in probability.
\item $\max |X_{nj}|\longrightarrow 0$ in probability.
\end{enumerate}
Then $\E(\exp(it S_n))\longrightarrow \exp(-\displaystyle\frac{t^2}2).$}\\

\begin{proof}
Let $t$ be fixed. From condition $(3)$ and $(4)$,
\[
|\sum_j r(tX_{nj})| \leq |t|^3 \sum_{j}|X_{nj}|^3 \\
\leq |t|^3 \max_j|X_{nj}|  \sum_j X_{nj}^2
\tend{n}{\infty}0 {\rm {~in~probability}}.
\]
\noindent{}$U_n=\exp\left (-\displaystyle \frac{t^2}2 \sum_j
X_{nj}^2+\sum_{j}r(tX_{nj})\right ) \longrightarrow  \exp(-\displaystyle
\frac{t^2}2)$ in probability as $n \longrightarrow +\infty$. This verifies the condition $(1)$ of Lemma
3.3 with $a=\exp(-\displaystyle \frac{t^2}2).$ It is easy to check
that the conditions $(2)$ ,$(3)$ and $(4)$ of the lemma 3.3 hold.
Thus $E\left( \exp(itS_n)\right) =E\left( T_nU_n\right) \longrightarrow
\exp(-\displaystyle \frac{t^2}2)$.
\end{proof}\\

\noindent{}Let $m$ be a positive integer and put

\begin{eqnarray*}
W_m \stackrel {\rm {def}}{=}&& \big \{ \big ({\sum_{i \in
I}\varepsilon_i p_i}\big ) h_m+\sum_{i \in I}\varepsilon_i
s_m(p_i) ~~:~~
\varepsilon_i \in \{0,-1,1\}, I \subset \{0,\cdots,p_m-1\},\\
&&p_i \in \{0,\cdots,p_m-1\} \}.
\end{eqnarray*}

\noindent{} The element $w =(\sum_{i \in
I}\varepsilon_ip_i)h_m+\sum_{i \in I}\varepsilon_is_m(p_i)$ is
called a word.\\

We shall need the following combinatorial lemma.\\

\noindent {\bf Lemma 3.5.}{ Under the
assumptions (i) and (ii) of theorem  2.1, all the words of $W_m$ are distinct.}\\

\begin{proof}{}
Let $w,w' \in W_m$, write
\begin{eqnarray*}
w &=&(\sum_{i \in I}\varepsilon_ip_i)h_m+\sum_{i \in
I}\varepsilon_is_m(p_i)\\
w' &=&(\sum_{i \in I'}\varepsilon'_ip'_i)h_m+\sum_{i \in
I}\varepsilon'_is_m(p'_i).
\end{eqnarray*}
\noindent{}Then $w=w'$ implies
\[
\{(\sum_{i \in I}\varepsilon_ip_i)-(\sum_{i \in
I'}\varepsilon'_ip'_i)\}h_m=\sum_{i \in
I'}\varepsilon'_is_m(p'_i)-\sum_{i \in I}\varepsilon_is_m(p_i)
\]
\noindent{}But the modulus of LHS is greater than $h_m$ and the
modulus of RHS is less than $2 \sum_{j=0}^{p_m-1}a_j^{(m)}$. It
follows that
\begin{eqnarray*}
\sum_{i \in I}\varepsilon_ip_i&=& (\sum_{i \in
I'}\varepsilon'_ip'_i){\rm {~~~and }} \\\sum_{i \in
I}\varepsilon_is_m(p_i)&=&\sum_{i \in I'}\varepsilon'_is_m(p'_i).
\end{eqnarray*}
Since from (i) we have  $s_m(p+1) \geq 3 s_m(p)$, we get that the
representation in the form $\sum_{i \in I}\varepsilon_is_m(p_i)$
is unique and the proof of the lemma is complete.
\end{proof}
\\

\begin{proof}{of proposition 3.2} Let $A$ a Borel set,
$\lambda(A)>0$. Using the Helly theorem we may assume that the
sequence $\displaystyle \frac{\sqrt{2}}{\sqrt{p_m}}
\sum_{j=0}^{p_m-1} \cos((jh_m+s_m(j))t)$ converge in distribution.
As is well-known, it is sufficient to show that for every real
number $x$,

\begin{eqnarray}{\label {eq:eq7}}
  \displaystyle \frac1{\lambda(A)}\int_A \exp\left \{-ix\frac{\sqrt{2}}{\sqrt{p_m}}
  \sum_{j=0}^{p_m-1} \cos((jh_m+s_m(j))t) \right \} dt
  \tend{m}{\infty} \exp(-\frac{x^2}2). \nonumber
\end{eqnarray}

\noindent{} To this end we apply theorem 3.4. in the following
context. The measure space is the given Borel set $A$ of positive Lebesgue
measure in the circle with the normalised measure and the random
variables are given by
$$
X_{mj}=\frac{\sqrt{2}}{\sqrt{p_m}} \cos((jh_m+s_m(j))t),~~~~{\rm {where}}~~~0 \leq j \leq p_m-1,~
m \in \N.
$$

\noindent{}It is easy to check that the variables $\{X_{mj}\}$
satisfy conditions (1) and (4). Further, condition (3) follows from fact that
$$
\int_0^{2\pi} \left |\sum_{j=0}^{p_m-1} X_{mj}^2-1 \right |^2 dt
\tend{m}{\infty} 0.
$$
\noindent{} It remains to verify comdition (2) of theorem 3.4. it is sufficient to show that
\begin{eqnarray}{\label {eq8}}
\int_A \prod_{j=0}^{p_m-1}\left (
1-ix\frac{\sqrt{2}}{\sqrt{p_m}}\cos((jh_m+s_m(j))t)\right)dt
\tend{m}{\infty}\lambda(A).
\end{eqnarray}

\noindent{}Write
\begin{eqnarray*}
\Theta_m(x,t)&=&\prod_{j=0}^{p_m-1}
\left(1-ix\frac{\sqrt{2}}{\sqrt{p_m}}\cos((jh_m+s_m(j))t\right)\\
&=&1+\sum_{w=1}^{N_m}{\rho_w}^{(m)}(x) \cos(wt),
\end{eqnarray*}

\noindent{}where $\rho_w=0$ if $w$ is not of the form $(\sum_{i
\in I}\varepsilon_ip_i)h_m+\sum_{i \in I}\varepsilon_is_m(p_i)$,
$N_m=\displaystyle \left \{\frac{p_m(p_m-1)}2\right
\}h_m+s_m(p_m-1)+s_m(p_m-2)+\cdots+1$.\\

\noindent{}We claim that it is sufficient to prove the following:

\begin{eqnarray}{\label{density}}
\int_{0}^{2\pi} R(t) \prod_{j=0}^{p_m-1} \left
(1-ix\frac{\sqrt{2}}{\sqrt{p_m}}\cos((jh_m+s_m(j))t)\right) dt
\tend{m}{\infty}\int_{0}^{2\pi} R(t) dt,
\end{eqnarray}

\noindent{}for any trigonometric polynomial $R$. In fact, assume
that (\ref{density}) holds and let $\epsilon>0$.
Then, by the density of trigonometric polynomials, one can find a
trigonometric polynomial $R_\epsilon$ such that
\[
||\chi_A-R_\epsilon||_1 <\epsilon,
\] where
\noindent{}$\chi_A$  is indicator function of $A$. But

\begin{eqnarray*}
\left |\prod_{j=0}^{p_m-1}\left (
1-ix\frac{\sqrt{2}}{\sqrt{p_m}}\cos((jh_m+s_m(j))t)\right)\right|  \leq
\left \{\prod_{j=0}^{p_m-1}\left ( 1+\frac{2x^2}{p_m}\right) \right
\}^{\frac12},
\end{eqnarray*}

\noindent{}Since $1+u \leq e^u$, we get

\begin{eqnarray}{\label{eq:eq11}}
\mid\Theta_m(x,t)\mid \leq e^{x^2}.
\end{eqnarray}

\noindent{}Hence, according to (\ref{density}), for $m$
sufficiently large, we have

\begin{eqnarray*}{\label{eq :eq12}}
&&\left |\int_A \Theta_m(x,t) dt -\lambda(A)\right|
=|\int_A\Theta_m(x,t)dt -\int_{0}^{2\pi} \Theta_m(x,t) R_\epsilon(t) dt+
\\&&\int_{0}^{2\pi} \Theta_m(x,t) R_\epsilon(t) dt-\int_{0}^{2\pi}  R_\epsilon(t) dt + \int_{0}^{2\pi}
R_\epsilon(t) dt-\lambda(A)|< e^{x^2}\epsilon+2\epsilon.
\end{eqnarray*}

\noindent{}The proof of the claim is complete. It still remains to prove
(\ref{density}). Observe that

$$\int_{0}^{2\pi} \Theta_m(x,t)R(t) dt =
  \int_{0}^{2\pi} R(t) dt+\sum_{w=1}^{N_m}{\rho_w}^{(m)}(x) \int_{0}^{2\pi} R(t)
  \cos(wt)dt$$

\noindent{}and for
$w=p_{i_1}h_m+s_m(p_{i_1})+\sum_{j=1}^{r}\varepsilon_j
\{(p_{i_j}h_m)+s_m(p_{i_j})\}$, we have

\[
|{\rho_w}^{(m)}(x)| \leq \frac{2^{1-r}|x|^r}{p_m^{\frac{r}2}},
\]

\noindent{}hence
\[
\max_{w \in W}|{\rho_w}^{(m)}(x)| \leq
\frac{|x|}{p_m^{\frac{1}2}}\tend{m}{\infty}0.
\]
\noindent{}Since $\displaystyle \sum_{w \in\Z}|\int_{0}^{2\pi} e^{-iwt}
R(t)dt|$ is bounded, we deduce
\[
|\sum_{w=1}^{N_m}{\rho_w}^{(m)}(x) \int_{0}^{2\pi} R(t)
  \cos(wt)dt|\leq \frac{|x|}{p_m^{\frac{1}2}}
  \sum_{w \in\Z}|\int_{0}^{2\pi} e^{-iwt} R(t)dt|\tend{m}{\infty}0.
\]
\noindent{}The proof of the proposition 3.2. is complete.
\end{proof}\\

\begin{proof}{of proposition 3.1}
Let $A$ be a Borel subset of $\T$, and $x \in ]1,+\infty[$, then, for any
positive integer $m$, we have

\begin{eqnarray*}{\label{eq :eqfin}}
\int_A ||P_m(\theta)|-1| d\lambda(\theta) &\geq& \int_{\{\theta
\in A
~~~:~~~~|P_m(\theta)|>x\}}||P_m|-1| d\lambda(\theta)\\
&\geq&(x-1)\lambda\{\theta \in A~:~|P_m(\theta)|>x\}\\
&\geq& (x-1) \lambda\{\theta \in A ~~:~~|\Re({P_m(\theta)})|>x\}
\end{eqnarray*}
\noindent{} Let $m$ go to infinity and use  propositions
2.5 and 3.2 to get

\[
\int_A \phi d\lambda \geq (x-1)\{1-{\mathcal
{N}}([-\sqrt{2}x,\sqrt{2}x])\}\lambda(A).
\]

\noindent{}Put $K=(x-1)\{1-{\mathcal {N}}([-\sqrt{2}x,\sqrt{2}x])\}$. Hence $
\alpha(A) \geq K \lambda(A)$, for any Borel subset $A$ of $\T$ which proves the proposition.
\end{proof}\\

\noindent{}Now, We give the proof of our main result.\\

\begin{proof}{of theorem 2.1} Follows easily from the proposition 2.6 combined with
proposition 3.1. 
\end{proof}\\

\noindent{}Let us mention that the same proof works for the
following more general statement.\\

\noindent {\bf Theorem 3.6. }{\it Let $T=T_{(p_k,
(a_j^{(k)})_{j=1}^{p_k})_{k=0}^\infty}$ be a rank one
transformation such that,

\begin{eqnarray}{\label{eq:eqt2}}
&&(i)~~a_{j+1}^{(k)}\geq 2s_k(j),j=0,\cdots, p_k-1,\nonumber\\
&& {\rm {~~with~~~~}} s_k(j)=a_1^{(k)}+\ldots +a_j^{(k)},
s_k(0)=0,\nonumber\\
&&(ii)~~\frac{s_k{(p_k)}-p_k\min_{1 \leq j \leq
p_k}(a_j^{(k)})}{h_k+\min_{1 \leq j \leq  p_k}(a_j^{(k)})}<\frac12
\nonumber
\end{eqnarray} then the spectrum of $T$ is singular.}
\\

\noindent{}{\bf Remark.} We note that, in \cite{Bourgain} and
\cite{Klemes1}, the strategy of the authors is to show that the
absolutely continuous measure  $\beta$, obtained as the limit of
some subsequence of the sequence $(||P_m|^2-1| d\lambda)_{m \geq
0}$, is equivalent to Lebesgue measure, in fact the authors prove
that
$$ \beta \geq K \lambda,~~~~~~~~~~~~~~~~ (E)$$
\noindent{}for some $K>0$. In the case of Ornstein
transformations, the relation (E) hold almost surely.\\

\section*{4. Simple Proof Of Bourgain Theorem.}

Bourgain Theorem deals with Ornstein transformations. In
Ornstein's construction, the $p_k$'s are rapidly increasing, and
the number of spacers, $a_i^{(k)}$, $1 \leq i\leq p_k-1$, are
chosen randomly. This may be organized in different ways as
pointed out by J. Bourgain in \cite{Bourgain}. Here we suppose that we are given
two sequences $(t_k)$, $(p_k)$ of positive integers and a sequence ($\xi_k$) 
 of probability measure such that the support of each
$\xi_k$ is a subset of $X_k = \{-\displaystyle
\frac{t_k}{2},\cdots,\displaystyle \frac{t_k}{2}\}$. We choose now
independently, according to $\xi_k$, the numbers
$(x_{k,i})_{i=1}^{p_k-1}$, and $x_{k,p_k}$ is chosen
deterministically in $\N$. We put, for $1 \leq i \leq p_k$,

$$a_i^{(k)} = t_{k} + x_{k,i} - x_{k,i-1}, ~~{\rm with} ~~x_{k,0}
= 0.$$\\

\noindent{} We have

$$h_{k+1} = p_k(h_k + t_{k}) + x_{k,p_k}.$$

\noindent{}So the deterministic sequences of positive integers
$(p_k)_{k=0}^\infty$, $(t_k)_{k=0}^\infty $ and
$(x_{k,p_k})_{k=0}^\infty$ determine completely the sequence of
heights $(h_k)_{k=0}^\infty$. The total measure of the resulting
measure space is finite if
\begin{eqnarray}\label{eqn:fini}
\sum_{k=0}^{\infty}\frac{t_k}{h_k}+\sum_{k=0}^\infty
\frac{x_{k,p_k}}{p_kh_k} < \infty.
\end{eqnarray}
  \noindent{}We will assume that this
requirement is satisfied.\\
We thus have a probability space of Ornstein transformations
$\Omega=\prod_{l=0}^\infty X_l^{p_l-1}$ equipped with the natural
probability measure $\P \stackrel {\rm def}
{=}\otimes_{l=1}^{\infty} P_l$, where $P_l\stackrel {\rm def}
{=}\otimes_{j=1}^{p_l-1}{\xi_j}$; ${\xi_j}$ is the probability
measure on $X_j$. We denote this space by $(\Omega, {\mathcal
{A}}, {\P})$. So $x_{k,i}$, $1 \leq i \leq p_k -1$, is the
projection from $\Omega$ onto the $i^{th}$ co-ordinate space of
$\Omega_k \stackrel {\rm def} {=} X_k^{p_k-1}$, $1 \leq i \leq
p_k-1$. Naturally each point $\omega =(\omega_k =
(x_{k,i}(\omega))_{i=1}^{p_k-1})_{k=0}^\infty$ in $\Omega$ defines
the spacers and therefore a rank one transformation
$T_{\omega,x}$, where $x=(x_{k,p_k})$.

\noindent{}This construction is more general than the
 construction due to Ornstein \cite{Ornstein} which corresponds to the case 
$t_k=h_{k-1}$, $\xi_k$ is the uniform distribution on $X_k$ and  $p_k >> h_{k-1}$. \\

\noindent{} We recall that Ornstein in \cite{Ornstein} proved that
there exists a sequence ${(p_k,x_{k,p_k})}_{k \in \N}$ such that,
$T_{\omega,x}$ is almost surely mixing. Later in \cite{Prikhodko}
Prikhod'ko obtains the same result for some special choice of the
sequence of the distribution ${(\xi_m)}$ and recently, using the
idea of D. Creutz and C. E. Silva \cite{Creutz-silva} one can
extend this result to a large family of probability measures
associated to Ornstein construction. In our general construction,
according to (\ref{eqn:type1}) the spectral type of each
$T_{\omega}$ , up to a discrete measure, is given by
   \[
\sigma_{T_\omega }=\sigma^{(\omega)}
_{\chi_{B_0}}=\sigma^{(\omega)} =W ^{*}\lim \prod_{l=1}^N\frac
1{p_l}\left| \sum_{p=0}^{p_l-1}z^{p(h_l+t_l)+x_{l,p}}\right|
^2d\lambda.
\]

\noindent{} With the above notations, we state Bourgain theorem in
the following
form:\\

\bigskip
\noindent {\bf Theorem 4.1 (\cite{elabdal2}). }{\it  For every
choice of $(p_k), (t_k), (x_{k,p_k})$ and for any family of
probability measures ${\xi_k}$ on $X_k=\{-t_k,\cdots,t_k\}$ of
$\Z$, $ {k \in \N^*}$, for which
$$\sum_{s \in X_m}\xi(s)^2 \tend{m}{\infty}0,$$ \noindent{}the
associated generalized Ornstein transformations has almost surely
singular spectrum. i.e.,
\[
\P\{\omega~:~ \sigma^{(\omega)} ~\bot ~\lambda\}=1.
\]
\noindent{}where $\P \stackrel {\rm def} {=}\otimes_{l=0}^{\infty}
\otimes_{j=1}^{p_l-1}{\xi_l}$ is the probability measure on
$\Omega=\prod_{l=0}^\infty X_l^{p_l-1}$. } \vskip0.5cm

\noindent{}In the context of Ornstein's constrution, we state
proposition 2.5 in the following form:\\

\noindent {\bf Proposition 4.2. }{\it There exist a subsequence of
the sequence $\left (\left |\left| P_m(z)\right|-1 \right|\right)$
which converge weakly to some non-negative function
$\phi(\omega,\theta)$ which satisfies  $ \phi \leq $2, almost surely
with respect the $\P \bigotimes \lambda$.}
\\

\begin{proof}
Easy exercise.
\end{proof}\\

\noindent Put, for any $\omega\in \Omega$
\[
\alpha_{\omega}=\phi(\omega,\theta)d\lambda.
\]
We shall prove that $\alpha_{\omega}$ is equivalent to the
Lebesgue measure for almost all $\omega$. In fact, we have the
following proposition:
\\

\noindent {\bf Proposition 4.3. }{\it There exist an absolutely
positive constant $K$ such that
  for almost all $\omega$ we have
$$\alpha_{\omega} \geq K \lambda.$$}
\\

\noindent The proof is based on the following two lemmas.

\noindent {\bf Lemma 4.4. }{\it $\displaystyle lim \int \int
||P_m|-|P_m'|| d\theta d\P \tendn 0$, where $\displaystyle
P_m'(\theta)=P_m(\theta)-\int_{\Omega}P_m(\theta)d\P.$ }

\begin{proof} For any $m \in \N$, we have
\begin{eqnarray*}
\int \int ||P_m|-|P_m'|| d\theta d\P &\leq& \int \int |P_m-P_m'|
d\theta d\P \\ & = &
\int |\int P_m d\P| d\theta\\
&\leq& \int |\sum_{p=0}^{p_m-1}\frac1{\sqrt{p_m}} z^{p(h_m+t_m))}| |\sum_{s\in
X_{m}}\xi_m(s)z^s| dz.
\end{eqnarray*}
Hence by Cauchy-Schwarz inequality
\begin{eqnarray*}
\int \int ||P_m|-|P_m'|| d\theta d\P \leq \sum_{s \in X_m}\xi(s)^2
\tend{m}{\infty}0.
\end{eqnarray*}
The proof of the lemma is complete.
\end{proof}\\
\vskip0.5cm \noindent Now observe that we have
\[
\int_T |\sum_{s\in X_{m}}\xi_m(s)z^s|^2 dz=
\int_T |\sum_{s\in X_{m}}\xi_m(s)z^{2s}|^2 dz=\sum_{s \in X_m}\xi(s)^2
\tend{m}{\infty}0.
\]
So, we may extract a subsequence $(m_k)$ for which, for almost all
$t \in [0,2\pi)$, we have
\[
\sum_{s\in X_{m_k}}\xi_{m_k}(s)e^{ i s t} \tend{k}{\infty}0
{\rm {~~and~~}} \sum_{s\in X_{m_k}}\xi_{m_k}(s)e^{ i  s 2t} \tend{k}{\infty}0
\]
\noindent Define
\[
\Theta
\stackrel{def}{=}\{\theta ~:~ \sum_{s\in X_{m_k}}\xi_{m_k}(s)e^{i
 s \theta}\tend{k}{\infty}0 {\rm {~and~}} \sum_{s\in X_{m_k}}\xi_{m_k}(s)e^{ i  s 2t}
 \tend{k}{\infty}0\}
\]
\noindent Choose $m \in \{m_k\}$, $t \in \Theta$ and put, for $j \in
\{0,\cdots,p_m-1\}$,
\begin{eqnarray*}
Y_{m,j}(\omega)&=&\cos((j(h_m+t_m)+x_{m,j}(\omega))t)-\int
\cos((j(h_m+t_m)+x_{m,j}(\omega))t) d\P,\\
Z_{m,j}(\omega)&=&\sqrt{\frac{2}{p_m}}Y_{m,j}(\omega).
\end{eqnarray*}
\\
\noindent {\bf Lemma 4.5. }{\it For any fixed $t \in \Theta$, the
distribution of the sequence of random variables $\displaystyle
\sum_{j=0}^{p_m-1} Z_{m,j}(\omega)$ converges
to the Gauss distribution. That is,
\\
\begin{eqnarray}{\label{eq:beq6}}
\P \left \{\omega \in \Omega~~:~~
\sum_{j=0}^{p_m-1} Z_{m,j}(\omega) \leq x
\right \} \tend{m}{\infty}\frac1{\sqrt{2\pi}}\int_{-\infty}^{x}
e^{-\frac12u^2}du\stackrel{\rm {def}}{=}{\mathcal {N}}(]-\infty,x]).
\nonumber
\end{eqnarray}
}
\begin{proof}
Since the random variables are independent, centred and uniformly bounded by $\displaystyle
\frac{2\sqrt{2}}{\sqrt{p_m}}$, the conditions (1), (2) and (4) of the
theorem 3.4 are satisfied. We have also the following:
\begin{eqnarray*}
&&\E \left \{ \left(\sum_{j=1}^{p_m-1}{Z_{m,j}^2-1}\right) ^{2} \right\rbrace \\
&& = \frac{4(p_m-1)}{p_m^2}\E ( Y_{m,1}^4)+\frac{(p_m-1)(p_m-2)}{p_m^2}
(\E (2Y_{m,1}^2))^2-2\frac{p_m-1}{p_m} \E (2 Y_{m,1}^2)+1.
\end{eqnarray*}
\noindent{}But $\E(2Y_{m,1}^2) \tend{m}{\infty}1$, it follows
that the variables $\{Z_{m,j}\}$ satisfy  condition (3) of theorem 3.4.
Thus all the conditions of theorem 3.4 hold and we conclude that the distribution of
$\displaystyle \sum_{j=0}^{p_m-1} Z_{m,j}(\omega)$
converges to normal distribution.
\end{proof}\\

\begin{proof}{of proposition 4.3}
Let $A$ be a Borel subset of $\T$, $C$ a cylinder set in
$\Omega$, and $x \in ]1,+\infty[$, then, for any positive integer
$m$, we have

\begin{eqnarray*}{\label{eq :beqfin}}
&&\int_{A \times C} ||P_m(\theta)|-1| d\lambda(\theta) d\P \\
&\geq& \P(C) \int_{A \times \Omega} ||P'_m(\theta)|-1|
d\lambda(\theta) d\P
-\int||P_m|-|P_m'||d\P d\lambda\\
&\geq& \P(C) \int_{ \{|P'_m|>x\} \bigcap A \times
\Omega}||P'_m|-1| d\lambda(\theta)d\P
-\int||P_m|-|P_m'||d\P d\lambda\\
&\geq&\P(C) (x-1)\int_A\P\{|\Re
{(P'_m(\theta))}|>x\}d\lambda-\int||P_m|-|P_m'||d\P d\lambda
\end{eqnarray*}
\noindent{} Let $m$ go to infinity and combine lemmas 4.4 and
4.5 to get

\[
\int_{A \times C} \phi d\lambda d\P \geq (x-1)\{1-{\mathcal
{N}}([-\sqrt{2}x,\sqrt{2}x])\}\lambda(A) \P(C).
\]

\noindent{}Put $K=(x-1)\{1-{\mathcal {N}}([-\sqrt{2}x,\sqrt{2}x])\}$. Hence, for
almost all $\omega$, we have, for any Borel set $A \subset \T$,  $ \alpha_{\omega}(A) \geq K
\lambda(A)$, and the proof of the
proposition is complete.
\end{proof}\\

\begin{proof}{of theorem 4.1}
Follow easily from the proposition 4.3 combined with the
proposition 2.6.
\end{proof}
\\

\noindent{}{\bf Remark.} We point out that there exist a rank one 
mixing transformations on the space with finite measure satisfying
the condition of theorem 3.6, In fact, following the notations of section 4, one may define
the spacers in the Ornstein construction by
\[
 {a_j^{(k)}}={3^j} t_k+x_{k,j}-x_{k,j-1},
\]

\noindent{}and choose the sequence $(t_k)_{k \in \N}$ such that the measure of dynamical system is finite.
Thus the condition of theorem 3.6 hold and the class is
mixing almost surely.\\

\begin{center} {\bf Acknowledgements}
\end{center}
I would like to express my thanks to J-P. Thouvenot who posed to me
the problem of the singularity of the spectrum of the rank one
transformations, and also to B. Host, F. Parreau and F. Bassam for
many conversations on the subject I have had with them. My thanks
also to Andrew Granville and the organizers of SMS-NATO ASI 2005
Summer School.\\
I am also grateful to the referee for a number of valuable suggestions and remarks and for 
having correted the mistakes in the first version of this work.

\bibliographystyle{nyjplain}

\end{document}